\documentclass{article}
\usepackage{amssymb}

\usepackage{amsmath}
\usepackage{amsthm}

\begin{document}

\title {Nice Bounds for the Generalized Ballot Problem}
\author {Delong Meng}
\date{Massachusetts Institute of Technology\\ Email: delong13@mit.edu}
\maketitle

\begin{abstract}
This paper gives two sharp bounds for the generalized ballot problem with 
candidate $A$ receiving at least $\mu$ times as candidate $B$ for an 
arbitrary real number $\mu$.  
\end{abstract}

\paragraph*{Introduction} Suppose in an election candidate $A$
received $a$ votes and candidate $B$ received $b$ votes. We count
the votes one at a time in any of the $\binom {a+b}{a}$ possible
sequences. Let $a_r$ and $b_r$ denote the number of votes $A$ and
$B$ have after counting the $r^{th}$ vote where $1\le r\le a+b$
(notice that $a_r+b_r=r$). Let $\mu$ be any positive real number. We
call a sequence desirable if $a>\mu b$ and $a_r>\mu b_r$ for all
$r$. We call a sequence cute if $a\ge\mu b$ and $a_r\ge \mu b_r$ for
all $r$. Let $P$ denote the probability that a sequence is desirable
and $P^*$ denote the probability that a sequence is cute.

Several authors started several articles quoting several well-known
results of the Ballot Problem. For brevity, Andre [7] and Barbier[8]
discovered that if $\mu\in\mathbb N$, then $$P=\dfrac{a-\mu
b}{a+b}.$$ Aeppli[9] showed that if $\mu\in\mathbb N$, then
$$P^*=\dfrac{a-\mu b+1}{a+1}.$$ Finally in 1962
Takacs[4] took a giant leap and bravely proved that for an arbitrary
$\mu\in\mathbb R$,
$$P=\dfrac{a}{a+b}\sum_{j=0}^{b}C_j\dfrac{\dbinom b j}{\dbinom {a+b-1}{j}}$$
where $C_0=1$ and $C_j$ satisfies the following recurrence formula:
$$\sum_{j=0}^{k}C_j\dfrac{\dbinom k j}{\dbinom {\lfloor k\mu\rfloor+k-1}{j}}=0$$
for all positive integers $k$. This formula gives the exact value
for $P$. However we can hardly imagine how big this number really
is. Therefore this article proves the following bounds:
\begin{equation}\label{one}
\textsc{ Theorem 1. }{ \dfrac{a-\lfloor\mu b\rfloor}{a+b}\le P\le
\dfrac{a-\lfloor\mu\rfloor b}{a+b}}
\end{equation}
\begin{equation}\label{two}
\textsc{ Theorem 2. } { \dfrac{\lfloor a-\mu b+1\rfloor}{a+b}\le
P^*\le \dfrac{a+1-\mu b}{a+1}.}
\end{equation}

\noindent We prove the two upper bounds with the Pseudo-Reflection
Principle and the two lower bounds with Penetrating Analysis.

\paragraph*{Pseudo-Reflection Principle} Let's start with Theorem 1.
We look at the relationship between the undesirable sequence and the
sequence with $a$ votes for $A$ and $b-1$ votes for $B$. We call
this the Pseudo-Reflection Principle because the case $\mu=1$ is 
essentially the reflection principle.
\footnote{However we are not finding any bijections here. We are only 
counting the number of undesirable sequences.} Both Goulden[5] and
Renault[6] proved the equality case (when $\mu\in \mathbb N$). They
both considered the smallest $r$ such that $a_r\le \mu b_r$. However
this approach does not generalize to the case when
$\mu\in\mathbb R$.\\
Instead we consider the largest $r$ such that $a_r\le \mu b_r$. When
the $r^{th}$ vote is counted, we must have $a_r=\lfloor\mu
b_r\rfloor\le \lfloor\mu\rfloor b_r$. There are
$\binom{a_r+b_r}{a_r}$ such undesirable sequences. Now consider the
number of sequences with $a_r$ votes for $A$ but $b_r-1$ votes for
$B$. For each $r$ there are
$\binom{a_r+b_r-1}{b_r-1}$ such sequences.\\
Consider the operation of replacing the first $r$ votes in an
undesirable sequence with these $\binom{a_r+b_r-1}{b_r-1}$
sequences. This operation yields all sequences with $a$ votes for
$A$ but $b-1$ votes for $B$ because for any sequence with $a$ votes
for $A$ but $b-1$ votes for $B$, there must
exist an $r$ such that $a_r-1\le\mu (b_r+1)<a_r$.\\
Since
\begin{equation}\label{notouch}
\dbinom{a_r+b_r-1}{b_r-1}=\dfrac{b_r}{a_r+b_r}\dbinom{a_r+b_r}{a_r}
\ge \dfrac{1}{\lfloor\mu\rfloor+1}\dbinom{a_r+b_r}{a_r}
\end{equation}
we deduce that the number of undesirable sequence is
at most $\lfloor\mu\rfloor+1$ times the sequences with $a_r$ votes
for $A$ but $b_r-1$ votes for $B$. Therefore
$$P\cdot \dbinom{a+b}{a} \le
\dbinom{a+b}{a}-(\lfloor\mu\rfloor+1)\dbinom{a+b-1}{b-1}=\dfrac{a-\lfloor\mu\rfloor
b}{a+b}\dbinom{a+b}{a}.$$

\noindent\textit{Remark } The conditions for equality to hold are
not trivial. Dvoretzky [11] proved that that equality holds if and
only if $\mu$ is sufficiently close to $\frac a b$ or $\mu$ is
sufficiently close to an integer. See Dvoretzky [11] for the
precise definitions of sufficiently close.\\

Now we move on to Theorem 2. Notice that the upper bound for Theorem
2 is a trivial consequence of theorem 1 when $\mu$ is an integer.
(We can simply add one vote for $A$ in the beginning of the sequence
and then $a_r>b_r$.) But such a correlation does not give the sharp
bound in Theorem 2 for $\mu\in\mathbb R$.\\
Using the Pseudo-Reflection technique we similarly consider the
largest $r$ such that $a_r<\mu b_r$. We have $a_r=\lceil\mu
b_r\rceil-1$. This time, however, we compare the ugly (non-cute)
sequence to the sequence with $a+1$ votes for $A$ and $b-1$ votes
for $B$. We similarly replace the first $r$ votes with sequences of
$a_r+1$ votes for $A$ and $b_r-1$ votes for $B$. This operation
yields all possible sequences with $a_r+1$ votes for $A$ and $b_r-1$
votes for $B$.\\
Now we have
\begin{equation}\label{touch}
\dbinom{a_r+b_r}{b_r-1}=\dfrac{b_r}{a_r+1}\dbinom{a_r+b_r}{a_r} \le
\dfrac{1}{\mu}\dbinom{a_r+b_r}{a_r}
\end{equation}
which implies that the number of ugly sequence is at least $\mu$
times the number of sequences of $a_r+1$ votes for $A$ and $b_r-1$
votes for $B$. Therefore
$$P^*\cdot \dbinom{a+b}{a} \le
\dbinom{a+b}{a}-(\mu)\dbinom{a+b}{b-1}=\dfrac{a-\mu
b+1}{a+1}\dbinom{a+b}{a}.$$

\noindent \textit{Remark } We can translate these sequences into
lattice paths from the origin the point $(b, a)$. A desirable path
never touches the line $y=\mu x$, and a cute path never go below the
line. The inequality $(\ref{notouch})$ shows that the number of
undesirable paths is at least $\lfloor\mu\rfloor$ times the number
of paths from the point $(1,0)$ to $(b,a)$. The inequality
$(\ref{touch})$ shows that the number of ugly paths is at least
$\mu$ times the number of path from the point $(1,-1)$ to $(b,a)$.
But intuition does not help much in this problem since it involves
calculations and one to $\lfloor\mu\rfloor$ correspondence rather
than a pure bijection.\\

\paragraph*{Penetrating Analysis} We first prove the lower bound
for Theorem 2. We claim that at least $\lfloor a-\mu b\rfloor+1$ of
the $a+b$ cyclic permutations of any given sequence of votes are
cute. This method is called penetrating analysis in Mohanty [1].

For any given sequence, define the weighted partial sum as
$S_r=a_r-\mu b_r$. Note that the sequence is cute if and only if
$S_r\ge 0$ for all $r.$ Suppose $S_i$ is the minimum (if there are
multiple $i$s then we can take any of them). We cyclically permute
the first $i$ terms of the sequence to the end of the sequence. In
other words we erase the first $i$ terms and attach them to the end
of the sequence. Now let $S'$ be the weighted partial sum for this
new sequence. We finish the proof with three lemmas.\\

\textit{Lemma 1. This new sequence is cute.}\\

\textit{Proof: } If $r\le a+b-i$, then $S'_r=S_{a+b-i}-S_i\ge 0$. If
$r>a+b-i$, then $S'_r=S_{r-(a+b-i)}+S_{a+b}-S_i\ge 0$ because
$S_{a+b}\ge 0$. Therefore $S'_r\ge 0$ for all $r$.\\

\textit{Lemma 2. A cyclic permutation that begins with the $r^{th}$
term of this sequence is cute if $S'_r\le S'_t$ for all $r+1\le t
\le a+b.$ For convenience, we call such an $r$ and also the $r^{th}$
vote cute.}\\

\textit{Proof: } Let ${S''}$ denote the weighted partial sum for the
cyclic permutation that begins with the $r^{th}$ term of this
sequence. If $j\le a+b-r$, then ${S''_j}=S'_{a+b-j}-S'_r\ge 0$. If
$j>a+b-r$, then ${S''_j}=S'_{j-(a+b-r)}+S'_{a+b}-S'_r\ge 0$ because
$S'_r\le S'_{a+b}$ and $S'_{j-(a+b-r)}\ge 0$. Therefore ${S''_r}\ge
0$ for all $r$.\\

\textit{Lemma 3. There exist at least $\lfloor a-\mu b\rfloor+1$
cute votes.}\\

\textit{Proof: } Let $r_1<r_2<\ldots<r_k=a+b$ denote all the cute
votes. We have $S'_{r_k}=a-\mu b$ and $S'_{r_1}\le 1$. Since
$S'_{r+1}\le S'_{r}+1$, we must have
\begin{equation}\label{duh}
S'_{r_{i+1}}\le S'_{r_i+1}\le S'_{r_i}+1
\end{equation}
(If $S'_{r_{i+1}}> S'_{r_i+1}$, then there must exist another cute
vote between $r_i+1$ and $r_{i+1}$, which contradicts the definition
of $r_i$.) Therefore $k-1\ge S'_{r_k}-S'_{r_1}\ge a-\mu b-1$.
Because $k$ is an integer, we have two cases:
\begin{enumerate}
    \item If $a-\mu b-1$ is not an integer, then
    $k\ge\lceil a-\mu b\rceil=\lfloor a-\mu b+1\rfloor.$
    \item If $a-\mu b-1$ is an integer, then consider all $r$ such
    that $a_r-\mu b_r<0$. Since there are finitely many such
    negative values, there exist an $\epsilon$ such that
    $a_r-(\mu-\epsilon)b_r\ge 0$ implies $a_r-\mu b_r\ge 0$.
    Replacing $\mu$ with $\mu-\epsilon$ would not affect the number
    of cute sequences. Therefore
    $k\ge\lceil a-(\mu-\epsilon)b\rceil=a-\mu b+1.$
\end{enumerate}
Now we have shown that that at least $\lfloor a-\mu b+1\rfloor$ of
the $a+b$ cyclic permutations of any given sequence are cute.
Therefore $P^*\ge \dfrac{\lfloor a-\mu b+1\rfloor}{a+b}.$

For Theorem 1, we can similarly show that at least $a-\lfloor\mu
b\rfloor$ of the $a+b$ cyclic permutations of any given sequence are
desirable. Notice that $\lfloor a-\mu b+1\rfloor=a-\lfloor\mu
b\rfloor$ if $a-\mu b-1$ is not an integer, and $\lfloor a-\mu
b+1\rfloor=a-\lfloor\mu b\rfloor+1$ otherwise. We can imitate the
proof for Theorem 2 until the last step. If $a-\mu b-1$ is an
integer, there does not $\epsilon$ such that $a_r-(\mu-\epsilon)b_r>
0$ can guarantee that $a_r-\mu b_r>0$ because $a_r-\mu b_r$ can be
equal to $0$. Therefore we can only conclude that $P\ge
\dfrac{a-\lfloor\mu b\rfloor}{a+b}.$\\

\noindent \textit{Remark 1} We can also prove the lower bound by
induction on $b$. Again let's prove the result of Theorem 2, and we
can follow the same procedure for Theorem 1. Base case $b=1$ is
trivial. Suppose that the bound is valid for all positive integers
less than $b$. We first apply lemma 1 to permute any given sequence
into a cute sequence. Then we consider two situations:
\begin{enumerate}
\item If there exist a cute $r$ such that $0<r<b$, then we can cut the
sequence into two cute sequences. By the inductive hypothesis, the
number of cute votes no less than $\lfloor a_r-\mu
b_r+1\rfloor+\lfloor a-a_r-\mu (b-b_r)+1\rfloor\ge\lfloor a-\mu
b+1\rfloor.$
\item If there does not exist a cute $r$ such that $0<r<b$, then for
all cute $r$, we must have $b_r=0$ or $b_r=b$. The rest follows
trivially from $(\ref{duh}).$
\end{enumerate}

\noindent\textit{Remark 2} Notice that any cute sequence must start
with a vote for $A$. We can thus treat all the votes for $B$ in
between two votes for $A$ as a single block. Therefore the argument
is still valid even if $B$ receives any number of weighted votes as
long as the weights add up to $b$. We can reformulate the problem as
follow.
\begin {quote}
Suppose in an election $A$ received $a$ votes all weighted 1.
However $B$ received $b'$ weighted votes whose sum is $b$. Both $a$
and $b$ are integers, and $\mu$ is a real number. We count the
$a+b'$ votes in a random order. Let $a_r$ and $b_r$ denote the sum
of the weighted votes $A$ and $B$ have after counting the $r^{th}$
vote where $1\le r\le a+b'$. Define $P$ as the probability that
$a_r>\mu b_r$ for all $1\le r\le a+b'$ (desirable sequences). Then
$$\dfrac{a-\lfloor\mu b\rfloor}{a+b'}\le P \le \dfrac{a}{a+b'}$$
\end {quote}
The upper bound is because each desirable sequence must start with a
vote for $a$. This bound, although achievable, is extremely weak
compare to $(\ref{one})$ and $(\ref{two})$. Goulden [5] discussed
the equality case for the lower bound when $\mu=1$ and all the
weights are integers. In fact equality holds for the lower bound if
$\mu$ and all the weights are integers. Therefore if all the weights
are integers, then $P\le \dfrac{a-\lfloor\mu\rfloor b}{a+b'}$, but
we still cannot find a sharp upper bound for arbitrary weights.

\paragraph*{Further Thoughts} The inspiration of this paper is to
search for a closed formula for any $\mu\in\mathbb R$. Intuitively
such a formula probably doesn't exist. Even if it does we must use
more advanced techniques. We still have many unanswered questions in
this paper. For example, can we find a sharper bound for Theorem 1
and Theorem 2? Can we find an upper bound if given specific weights
in the last problem? Can we derive similar inequalities for a
multi-candidate election?

\subsection*{References}

1. Mohanty, Gopal, Lattice Path Counting and Applications. Acedamic
Press, New York, 1980.\\
2. Marc Renault, Four Proofs of the Bollot Theorem,
\textit{Mathematics Magazine} 80, December 2007, 345-352.\\
3. L. Takacs, Ballot Problems. Z. Wahrscheinlichkeistheorie, 1
(1962), 154-158. \\
4. L. Takacs, A generalization of the ballot problem and its
application to the theory of queues. J. Amer. Statist. Assoc.,
57(1962), 327-337\\
5. I.P. Goulden and Luis G. Serrano, Maintaining the Spirit of the
Reflection Principle when the Boundary has Arbitrary Integer Slope,
J. Combinatorial Theory (A) 104 (2003) 317-326.\\
6. Marc Renault, Lost (and Found) in Translation: Andre's Actual
Method and its Application to the Generalized Ballot Problem,
\textit{American Mathematical Monthly} 115 (2008) 358-363.\\
7. D. Andre, Solution directe du probleme resolution par M.
Bertrand, \textit{Comptes Rendus de l'Academie des Sciences}, Paris,
105
(1887) 436¨C437.\\
8. E. Barbier, Generalisation du probleme resolution par M. J.
Bertrand, \textit{Comptes Rendus de l'Academie des Sciences, Paris},
105
(1887) p. 407.\\
9. A. Aeppli, Zur Theorie Verketteter Wahrscheinlichkeiten. These,
Zurich (1924).\\
10. J. Saran and K. Sen, Some Distribution Results on Generalized
Ballot Problems. Svazek 30 (1985) 157-165.\\
11. A. Dvoretzky and Th. Motzkin, A Problem of Arrangements,
\textit{Duke Math. J.} \textbf{14} (1947) 305-313.

\end{document}